# EQUIVARIANT K-THEORY OF REAL VECTOR SPACES AND REAL PROJECTIVE SPACES

by Max Karoubi

The computation of the equivariant K-theory $K_G^*(V)$ of the Thom space of a <u>real</u> vector bundle has been done successfully only under some spinoriality conditions [1], thanks to a clever use of the Atiyah-Singer index theorem (even if G is a finite group). One purpose of this paper is to fill this gap, at least for real vector spaces (considered as vector bundles over a point). For this purpose, we use at the same time the results in [7] (generalizing those of Atiyah) and the equivariant Chern character (Slominska [10], Baum and Connes [4]).

The interest of such a computation comes from many sources. First, it answers a question raised recently by Pierre-Yves Le Gall and Bertrand Monthubert [9] in their investigations on a suitable index theorem for manifolds with corners. They have to compute an "indicial K-theory" which is precisely the equivariant K-theory $K_G^*(V)$, where $V = \mathbf{R}^n$ and $G = \mathfrak{S}_n$ the symmetric group of n letters acting naturally by permutation of the coordinates in $\mathbf{R}^n$.

Secondly, these topological computations are linked with two natural algebraic questions. One of them is the determination of the number of conjugacy classes of a subgroup of the orthogonal group O(V) which split in the central extension induced by the pinorial group Pin(V) (cf. 2.4). The second question is how to compute the number of simple factors in the crossed product algebra $G \ltimes C(V)$, where C(V) is the Clifford algebra of V (cf. 2.11, 3.6 and 3.11).

Finally, these methods enable us to determine completely the rank of the group $K_G^*(P(V))$, where P(V) is the <u>real</u> projective space of V, in terms of the number of certain conjugacy classes of G (cf. 3.4 and 3.8). As we shall see, Algebra and Topology are deeply linked in these computations.

The previous results have a pleasant formulation when G is the symmetric group $\mathfrak{S}_n$ acting on $V = \mathbf{R}^n$ as above. In this case (cf. 1.9), the groups $K_G^0(V)$ and $K_G^1(V)$ are <u>free</u> of rank[1] equal to $p_n$ and $i_n$ respectively. Here $p_n$ (resp. $i_n$) denotes the number of partitions of n of the type

$$n = \lambda_1 + ... + \lambda_{2k}$$

---

[1] This rank was already determined in the paper by Le Gall and Monthubert [9]. What we essentially prove here is the freeness of these two groups.

with $1 \leq \lambda_1 < ... < \lambda_{2k}$ (resp. of the type
$$n = \lambda_1 + ... + \lambda_{2k+1}$$
with $1 \leq \lambda_1 < ... < \lambda_{2k+1}$).

Following the above program, we are able to determine completely the number[2] of simple factors in the crossed product algebra $\mathfrak{S}_n \ltimes C(\mathbf{R}^n)$. We also determine the rank of the groups $K_G^1(P(V))$ and $K_G^0(P(V))$ : they are $p_n$ and $2P(n) - p_n - i_n$ respectively, where $P(n)$ denotes the total number of partitions of n.

As an irony of mathematical history, the Euler-Poincare characteristic $p_n - i_n$ associated to this equivariant K-theory of $\mathbf{R}^n$ happens to have been determined by... Euler. It may be deduced from his well known "pentagonal identity" [3]

$$\sum_{n=1}^{\infty} (p_n - i_n) q^n = \prod_{m=1}^{\infty} (1 - q^m) = \sum_{m=-\infty}^{+\infty} (-1)^m q^{\frac{m(3m-1)}{2}}$$

In particular, $p_n - i_n$ is always equal to 0, 1 or -1. On the other hand, the integers $p_n$ and $i_n$ increase very fast with n and satisfy the following asymptotic expansion :

$$p_n \sim i_n \sim \frac{e^{\pi\sqrt{\frac{n}{3}}}}{8.3^{\frac{1}{4}}.n^{\frac{3}{4}}}$$

What we have just said is included in the first three sections of the paper. In the fourth section, we work out the R(G)-module structure of $K_G^*(V)$, when $G = \mathfrak{S}_n$ and $V = \mathbf{R}^n$. This structure is much simpler if we work over $\mathbf{Q}$, thanks to the Adams operations. When G is an arbitrary finite group, we may also determine $K_G^*(V) \otimes \mathbf{C}$ as a $R(G) \otimes \mathbf{C}$-module, thanks to the results developed by Slominska [10] and Baum-Connes [4], quoted at the beginning.

---

[2] This number is $p_n + 2i_n$.

[3] See for instance the book of G.E. Andrews : The theory of partitions (Addison-Wesley, 1976), p. 11 and also p. 97 for the asymptotic expansion considered in the same introduction.

# 1. Equivariant K-theory of real vector spaces.

**1.1.** Let X be a compact space on which a finite group G acts and let V be a <u>real</u> vector bundle over X provided with a linear action of G, compatible with the projection on X. We can define on V a <u>positive</u> definite metric invariant by the action of G and consider the associated Clifford bundle C(V) ; the group G also acts naturally on C(V). We denote by $\mathcal{E}_G^V(X)$ the category of real vector bundles where G and C(V) act simultaneously ; these two actions are linked together by the formula

$$g*(a.e) = (g*a).(g*e)$$

where the symbol $*$ (resp. $.$) denotes the action of G (resp. of C(V)). If X is reduced to a point, this category is simply denoted by $\mathcal{E}_G^V$.

**1.2.** The method developped in [8] § 1 for instance shows that $\mathcal{E}_G^V(X)$ is equivalent to the category of finitely generated projective modules over the crossed product algebra $G \ltimes \tilde{C}(V)$, where $\tilde{C}(V)$ is the algebra of continuous sections of the bundle C(V). On the other hand, if "1" denotes the trivial bundle of rank one (with the trivial action of G), we have a "restriction" functor

$$\varphi : \mathcal{E}_G^{V \oplus 1}(X) \longrightarrow \mathcal{E}_G^V(X)$$

The following theorem, where V is considered as a locally compact space, is proved in [7] for both real and complex K-theory with compact supports :

**1.3. THEOREM.** *The group* $K_G^*(V)$ *is naturally isomorphic to the Grothendieck group* $K^*(\varphi)$ *of the Banach functor* $\varphi$ [6]. *If* $\mathcal{C}'$ *denotes the category* $\mathcal{E}_G^{V \oplus 1}(X)$ *and* $\mathcal{C}$ *the category* $\mathcal{E}_G^V(X)$, *we therefore have the exact sequence*

$$K^{i-1}(\mathcal{C}') \longrightarrow K^{i-1}(\mathcal{C}) \longrightarrow K_G^i(V) \longrightarrow K^i(\mathcal{C}') \longrightarrow K^i(\mathcal{C})$$

**1.4.** Let us suppose now that X is reduced to a point. The categories $\mathcal{C}$ and $\mathcal{C}'$ are then semi-simple (all exact sequences are split). Therefore, they are both equivalent to a category of finitely generated modules over a semi-simple algebra, i.e. a product of matrix algebras over **C** (in the complex case), over **C**, **R** or **H** (in the real case). In the complex case, it follows that $K^0(\mathcal{C}')$ is free of finite type and $K^1(\mathcal{C}) = 0$. Therefore, $K_G^0(V) = K^0(\varphi) = \text{Ker}[K^0(\mathcal{C}') \longrightarrow K^0(\mathcal{C})]$ is free of finite type. More generally, $K_G^{-p}(V)$ is the Grothendieck group of the functor

$$\mathcal{E}_G^{V \oplus p \oplus 1}(X) \longrightarrow \mathcal{E}_G^{V \oplus p}(X)$$

where "p" denotes the trivial bundle of rank p, with the trivial action of G. Therefore, the group $K_G^1(V) = K_G^{-p}(V)$ for p odd (in the complex case) is also free by the same argument.

**1.5. Remark.** Let us suppose that G acts on V via oriented automorphisms and that the rank n of V is even, say n = 2r. If $e_1, ..., e_n$ is an orthonormal basis of V, the product $(i)^r e_1 ... e_n$ in the Clifford algebra C(V) has square +1 and anticommutes with each $e_\alpha$. It follows that the category $\mathcal{E}_G^{V \oplus 1}(X)$ splits as the product $\mathcal{E}_G^V(X) \times \mathcal{E}_G^V(X)$, the functor $\varphi : \mathcal{E}_G^{V \oplus 1}(X) \longrightarrow \mathcal{E}_G^V(X)$ being identified with the "sum" functor. Therefore, the group $K_G^i(V)$ is isomorphic to the group $K^i$ of the Banach category $\mathcal{E}_G^V(X)$. In particular, if X is reduced to a point, $K_G^1(V) = 0$, in contrast with the corollary 1.9 below. In real K-theory (with X reduced to a point again), the same argument shows that $K_G^i(V) = 0$ if $i \equiv 1$ mod. 4 and is free of finite type if $i \equiv 0$ mod. 4.

**1.6.** It remains to compute more precisely the rank of the free groups $K_G^i(V)$ (we consider now only <u>complex</u> K-theory), for an arbitrary finite group G. In this direction, we may use the isomorphism shown by Slominska, Baum and Connes [10] [4][9] between the equivariant K-theory tensored by **C** and equivariant cohomology. More precisely, let <G> be the set of conjugacy classes of elements of G and $g_1, ..., g_p$ be a set of representatives. The vector space $K_G^*(V) \otimes \mathbf{C}$ is then isomorphic to the following direct sum of cohomology vector spaces

$$K_G^*(V) \otimes \mathbf{C} = \bigoplus_{g_i \in <G>} H_c^*(V^{g_i}/C_{g_i}; \mathbf{C})$$

where $C_{g_i}$ denotes the centralizer of $g_i$ and * denotes the degrees mod. 2.

**1.7.** Let us set the following definitions : a conjugacy class $<g_i>$ is called <u>even</u> (resp. <u>odd</u>) if the dimension of $V^{g_i}$ is even[4] (resp. odd). It is called <u>oriented</u>[5] if all the elements of $C_{g_i}$ act on $V^{g_i}$ by oriented automorphisms.

**1.8. THEOREM.** *Let* G *be an arbitrary finite group acting linearly on a <u>real</u> vector space* V *of finite dimension. The group* $K_G^0(V)$ *(resp.* $K_G^1(V)$*) is then a free* **Z***-module of rank the number of conjugacy classes* $<g_i>$ *which are oriented and even (resp. odd).*

---

[4] The dimension 0 is not excluded.

[5] We may also say <u>positively oriented</u>, if we wish to follow the conventions used further in the paper (cf. 3.3).

*Proof.* The freeness of the groups has already been shown in 1.4. On the other hand, the $H_c^*(V^{g_i}/C_{g_i} ; \mathbf{C})$ are the <u>reduced</u> cohomology vector spaces in degrees $q = \dim(V^{g_i})$ of $S^q/\Gamma$, where $\Gamma = C_{g_i}$. It is easy to see that each vector space is isomorphic to the invariant part of $\tilde{H}^q(S^q ; \mathbf{C})$ under the action of $\Gamma$. Therefore, its dimension is 1 if $g_i$ is an oriented conjugacy class and 0 otherwise.

The following corollary (where $i_n$ and $p_n$ are defined in the introduction) is an extension of a theorem of Le Gall and Monthubert [9]. It may be deduced from the preceding theorem if we decompose the elements $g_i$ of $\mathfrak{S}_n$ into a product of cycles : the oriented conjugacy classes are then in bijective correspondence with the products of cycles of different lengths.

**1.9. COROLLARY.** *Let* $G = \mathfrak{S}_n$ *acts naturally on* $V = \mathbf{R}^n$. *Then we have*
$$K_G^0(\mathbf{R}^n) \cong \mathbf{Z}^{p_n} \text{ and } K_G^1(\mathbf{R}^n) \cong \mathbf{Z}^{i_n}$$

**1.10.** Let us go back to the notations of 1.1. If T is a spinorial representation space of G of <u>even</u> dimension 2r, the categories $\mathcal{E}_G^W$ and $\mathcal{E}_G^{T \oplus W}$ are then equivalent. As a matter of fact, let M be an irreducible C(T)-module and let $\bar{\sigma} : G \longrightarrow \mathrm{Spin}(T)$ be a lifting of the homomorphism $\sigma : G \longrightarrow O(T)$, where $O(T)$ denotes the orthogonal group of T. The category equivalence $\phi : \mathcal{E}_G^W \longrightarrow \mathcal{E}_G^{T \oplus W}$ is defined by associating to an object E of $\mathcal{E}_G^W$ the object $M \otimes E$ of $\mathcal{E}_G^{T \oplus W}$. In this formula, G acts on M via the representation $\bar{\sigma}$ ; the action of the Clifford algebra $C(T \oplus W)$ is induced from the homomorphism $T \oplus W \longrightarrow \mathrm{End}(M \otimes E)$ defined by
$$(t, w) \mapsto t \otimes 1 + \varepsilon \otimes w$$
where $\varepsilon = i^r e_1 \ldots e_{2r}$ and where $(e_\alpha)$ is an <u>oriented</u> orthonormal basis of T. We remark that this equivalence $\phi$ may be extended to the associated $\mathbf{Z}/2$-graded categories. Therefore, we also have the isomorphisms $K_G^0(W) \cong K_G^0(T \oplus W)$ and $K_G^1(W) \cong K_G^1(T \oplus W)$, according to 1.3. Here is an algebraic consequence :

**1.11. THEOREM.** *For any real G-vector space* W, *let us denote by* $a_0(W)$ *(resp.* $a_1(W)$*) the number of conjugacy classes of G which are even (resp. odd) with the notations of 1.7. Let us now consider two G-modules V and V' of the same dimension such that* $w_1(V) = w_1(V')$ *and* $w_2(V) = w_2(V')$, *where* $w_1$ *and* $w_2$ *denote the first two Stiefel-Whitney classes*[6]. *Then we have* $a_0(V) = a_0(V')$ *and* $a_1(V) = a_1(V')$. *In particular, if* $w_1(V') = w_2(V') = 0$, $a_i(V) = |<G>|$ *and* $a_{i+1}(V) = 0$ *where* $i = \dim V \bmod. 2$ *and where* $<G>$ *denotes the number of conjugacy classes of G.*

---
[6] Belonging therefore to $H^i(G ; \mathbf{Z}/2)$ with $i = 1, 2$.

*Proof.* According to 1.10, it is enough to find two spinorial representation spaces T and T' of even dimension such that V ⊕ T is isomorphic to V'⊕ T'. We choose

$$T = V \oplus V \oplus V' \oplus V' \text{ and } T' = V \oplus V \oplus V \oplus V'$$

The vector spaces T and T' are spinorial, since their first two Stiefel-Whitney classes are reduced to 0.

**1.12. Remark.** The R(G)-module structure of $K_G^*(V)$ is determined in Section 4 of the paper.

## 2. Relation with the conjugacy classes of the Schur group

**2.1.** Let us consider the 2-fold covering[7] Pin(n) of O(n), and the associated covering $\tilde{G}$ of G ⊂ O(n), via the pull-back diagram :

$$\begin{array}{ccc} \tilde{G} & \longrightarrow & \text{Pin}(n) \\ \downarrow & & \downarrow \\ G & \longrightarrow & O(n) \end{array}$$

The kernel of the homomorphism $\tilde{G} \longrightarrow G$ is identified with the multiplicative group ± 1 in the Clifford algebra $C(\mathbf{R}^n)$. If G is the symmetric group $\mathfrak{S}_n$, $\tilde{G}$ is the classical Schur group[8] which is well known in representation theory [5]. To this extension of G we associate a cocycle which is defined as follows. Let g and h be two elements of G which commute and let $\tilde{g}$, $\tilde{h}$ be two elements of $\tilde{G}$ above them. The commutator $[\tilde{g}, \tilde{h}] = \pm 1$ is not trivial in general. More generally, if $k_\alpha$ and $h_\alpha$ are finite families of elements of G such that the product of the commutators $[k_\alpha, h_\alpha]$ is equal to 1, the sign of the product $[\tilde{k}_\alpha, \tilde{h}_\alpha]$ defines the required cocycle. We would like to relate this sign to the "orientation" of the elements $g = g_i$ defined in 1.7.

**2.2.** Let $V_\lambda$ be the eigenspaces (of dimension 1 or 2) of g acting on V and let $h \in C_g$. We may put aside the non real eigenvalues, since we see by deformation of these eigenvalues to 1 ($V_\lambda$ being stable by h) that the commutator $[\tilde{g}, \tilde{h}]$ only depends on the eigenspaces of g and h associated to the eigenvalues +1 and -1. Let $V^+$ (resp. $V^-$) be the eigenspace of g associated to the eigenvalue +1 (resp. -1) which is of dimension $n^+$ (resp. $n^-$). We may notice that the element of the orthogonal group associated to the product $e_1 \ldots e_p$ in Pin(n), p ≤ n, is the

---

[7] In contrast with the conventions in [3], we consider here the Pin group associated to the Clifford algebra of $\mathbf{R}^n$ provided with the standard definite <u>positive</u> quadratic form (see also [6]).

[8] We keep this terminology for an arbitrary finite subgroup of the orthogonal group.

symmetry $e_i \mapsto \pm e_i$, where the sign is +1 if $i > p$ and -1 if $i \leq p$, a fact which enables us to compute the commutator $[\tilde{g}, \tilde{h}]$ easily. The argument of deformation and "crossing" of the eigenvalues quoted above reduces the discussion to one of the following four cases :

| det(g) | det(h) | $[\tilde{g}, \tilde{h}]$ | $\det(h|_{V^+})$ |
|---|---|---|---|
| $(-1)^{n^-}$ | +1 | +1 | +1 |
| $(-1)^{n^-}$ | -1 | $(-1)^{n^-}$ | -1 |
| $(-1)^{n^-}$ | -1 | $(-1)^{n^-+1}$ | +1 |
| $(-1)^{n^-}$ | +1 | -1 | -1 |

**2.3.** From this table we can decide easily for a given g whether there exists, or not, an element $h \in C_g$ such that the commutator $[\tilde{g}, \tilde{h}]$ is equal to -1 (which means that $\tilde{g}$ is conjugate to $-\tilde{g}$). Following the terminology of [5], we say that a conjugacy class <g> of G is decomposed in $\tilde{G}$ if $\tilde{g}$ is <u>not</u> conjugate to $-\tilde{g}$ in the Schur group $\tilde{G}$. The following theorem is then a consequence of the previous worksheet ; it is proved in [5] p. 29 when $G = \mathfrak{S}_n$. We here use the following notations : $h^+ = h|_{V^+}$, with $V^+ = \mathrm{Ker}(\rho(g) - 1)$.

**2.4. THEOREM.** *The conjugacy classes <g> of G decomposed in $\tilde{G}$ may be divided into two categories*

   a) *The permutation g is even. Then $\forall\, h \in C_g$, we have $\det(h) = \det(h^+)$. In particular, if $G = \mathfrak{S}_n$, g may be decomposed into a product of cycles of odd lengths.*

   b) *The permutation g is odd. Then $\forall\, h \in C_g$, we have $\det(h^+) = 1$. In particular, if $G = \mathfrak{S}_n$, g is decomposed into a product of cycles of different lengths.*

**2.5. DEFINITION-THEOREM.** *Let $\sigma : \tilde{G} \longrightarrow GL(W)$ be a representation of the group $\tilde{G}$. We say that $\sigma$ is of linear type if it satisfies the additional assumption $\sigma(-\tilde{g}) = -\sigma(\tilde{g})$. The (finite) set of isomorphism classes of irreducible representations of linear type is then in bijective correspondence with the conjugacy classes of elements in G which are decomposed in $\tilde{G}$.*

*Proof.* Since -1 is in the centre of $\tilde{G}$, it is clear that the irreducible representations of $\tilde{G}$ are divided into two categories : those of linear type and those which satisfy the property $\sigma(-\tilde{g}) = \sigma(\tilde{g})$ for all $\tilde{g}$, i.e. coming from representations of G. According to character theory, the first ones are in duality with the conjugacy classes of elements in G which are not decomposed in $\tilde{G}$, while the second ones are in duality with two copies of the set of conjugacy classes of G which are decomposed in $\tilde{G}$. Therefore, we have the identity

$$|<\tilde{G}>| = |<G>| + |<G>_{déc.}|$$

where $<G>_{déc.}$ denotes the set of conjugacy classes of elements in G which are decomposed in $\tilde{G}$.

**2.6. THEOREM.** *Let* V *be a vector space of even dimension provided with a positive definite quadratic form. Then the rank of the free group*[9] $K(\mathcal{E}_G^V)$ *is the number of conjugacy classes of elements in* G *which are decomposed in* $\tilde{G}$.

*Proof.* Let $V = \mathbf{R}^{2r}$, $\varepsilon = i^{(r+1)} e_1 \ldots e_{2r}$ and j be the canonical inclusion of Pin(V) in C(V)*, the group of invertible elements in the Clifford algebra C(V). We define <u>another</u> homomorphism $\phi$ from Pin (V) to C(V)* by putting $\phi(u) = j(u)$ if u is of even degree and $\phi(u) = j(u).\varepsilon$ if u is of odd degree. Since $(\varepsilon)^2 = -1$ and $\varepsilon$ anticommutes with the generators of the Clifford algebra, $\phi$ is a group homomorphism. On the other hand, we have a commutative diagram

$$\begin{array}{ccc} \text{Pin}(V) & \xrightarrow{\phi} & C(V)^* \\ \downarrow \theta & & \downarrow \Theta \\ O(V) & \longrightarrow & \text{Aut}(C(V)) \end{array}$$

In this diagram, $\theta$ is the standard covering of the orthogonal group by the Pin group and $\Theta$ associates to $u \in C(V)^*$ the inner automorphism $s \mapsto u.s.u^{-1}$. Therefore, we also have the commutative diagram

$$\begin{array}{ccc} \tilde{G} & \xrightarrow{\phi} & C(V)^* \\ \downarrow & & \downarrow \Theta \\ G & \longrightarrow & \text{Aut}(C(V)) \end{array}$$

which we shall use to "untwist" the action of G on C(V). More precisely, we are going to show that the category $\mathcal{E}_G^V$ is equivalent to the category $\mathcal{E}_{\tilde{G}}^{(V)}$ of finite vector spaces with an action of $\tilde{G}$ of <u>linear type</u> and an action of the Clifford algebra C(V) commuting with it. Since C(V) is isomorphic to a matrix algebra, it will follow (by Morita equivalence) that the category $\mathcal{E}_G^V$ is equivalent to the category of finite vector spaces provided with an action of $\tilde{G}$ which is of linear type.

---

[9] The definition of the category $\mathcal{E}_G^V$ is given in 1.1. Note that the rank of $K(\mathcal{E}_G^V)$ is also the number of simple factors in the crossed product algebra $G \ltimes C(V)$ which is semi-simple.

In order to carry out this program, let us consider an object M of $\mathcal{E}_G^V$, with an action of G denoted by $(g, m) \mapsto g*m$ and an action of C(V) denoted by $(\lambda, m) \mapsto \lambda.m$, as in 1.1. We may associate to M an object of the category $\mathcal{E}_{\tilde{G}}^{(V)}$ by keeping the same action of the Clifford algebra, but with a new action $\odot$ of $\tilde{G}$, defined by the following formula

$$\tilde{g} \odot m = \phi(\tilde{g})^{-1}.(g*m) \qquad (F)$$

where g denotes the projection of $\tilde{g}$ on G. According to the definitions in 1.1, we have the identity $g*(\phi(\tilde{g})^{-1}.m) = \phi(\tilde{g})^{-1}.(g*m)$. On the other hand, if $\lambda$ is an element of the Clifford algebra, we have $\tilde{g} \odot (\lambda.m) = g*(\phi(\tilde{g})^{-1}\lambda m) = \phi(\tilde{g})(\phi(\tilde{g})^{-1}\lambda) \phi(\tilde{g})^{-1} g*m = \lambda.(\tilde{g} \odot m)$, which shows that the actions of $\tilde{G}$ and C(V) commute. Finally, if -1 is the non trivial central element of $\tilde{G}$ above 1 in G, we have $(-1) \odot m = \phi(-1)^{-1}.(1*m) = -m$, which shows that the action of $\tilde{G}$ on M is of linear type. The previous considerations enable us to define a functor from $\mathcal{E}_G^V$ to $\mathcal{E}_{\tilde{G}}^{(V)}$. The functor in the other direction is defined in the same way. This concludes the proof of the theorem.

**2.7.** The previous proof may be adapted to the category $\mathcal{E}_G^{V \oplus 1}$, if V is even dimensional as before, regarding the objects of $\mathcal{E}_G^{V \oplus 1}$ as <u>graded</u> objects of the category $\mathcal{E}_G^V$. More precisely, with the notations as in 2.6, we give ourselves an involution $\eta$ of the object M such that $\eta.v = -v.\eta$ for $v \in V$ and $g.\eta = \eta. g$ for $g \in G$. According to the formula (F) above, the equivalence betweeen the categories $\mathcal{E}_G^V$ et $\mathcal{E}_{\tilde{G}}^{(V)}$ changes the action of G into an action of $\tilde{G}$ such that $\eta.\tilde{g} = \tilde{g}.\eta$ if $\det(g) = 1$ and $\eta.\tilde{g} = - \tilde{g}.\eta$ if $\det(g) = -1$. On the other hand, by changing $\eta$ into $\eta.\varepsilon$, we note that $\eta.\varepsilon$ commutes with the action of C(V) and is still an involution. In other terms, the category $\mathcal{E}_G^{V \oplus 1}$ is equivalent to the category of finite vector spaces E with an action of $\mathbf{Z}/2 \ltimes \tilde{G}$ of linear type, if we view this last group as the "Schur group" $\tilde{G}'$ associated to the subgroup $G' = \mathbf{Z}/2 \times G$ of $O(n+1)$.

Finally, we notice that the case where V is odd dimensional follows from the even case, thanks to the category equivalence $\mathcal{E}_G^V \approx \mathcal{E}_G^{(V \oplus 1) \oplus 1}$ described in 1.10. Hence, we have the following theorem :

**2.8. THEOREM.** *Let* G *be a subgroup of* $O(n)$ *and let* $V = \mathbf{R}^n$. *The ranks* $R_G^V$ *and* $R_G^{V \oplus 1}$ *of the groups* $K(\mathcal{E}_G^V)$ *and* $K(\mathcal{E}_G^{V \oplus 1})$ *respectively are determined as the following functions of* n :

1. *If* n *is even (resp. odd),* $R_G^V$ *(resp.* $R_G^{V\oplus 1}$*) is the number of conjugacy classes of* G *which are decomposed in* $\tilde{G}$

2. *If* n *est odd (resp. even),* $R_G^V$ *(resp.* $R_G^{V\oplus 1}$*) is the number of conjugacy classes of the subgroup* G' = **Z**/2 x G *of* O(n+1) *which are decomposed in* $\tilde{G}'$ =**Z**/2 ⋉ $\tilde{G}$.

**2.9. COROLLARY** *(compare with 1.11).* *Let* ρ : G ⟶ O(n) *and* ρ' : G ⟶ O(n) *be two representations of* G *into the orthogonal group* O(n). *Let* $w_i$ *(resp.* $w'_i$*) be the Stiefel-Whitney classes of* ρ *and* ρ' *respectively. We assume that* $w_1 = w'_1$ *and* $w_2 = w'_2$. *Finally, let* $\tilde{G}$ *and* $\tilde{G}'$ *be the associated Schur groups. Then the number of conjugacy classes of* G *which are decomposed in* $\tilde{G}$ *is equal to the number of conjugacy classes of* G' *which are decomposed in* $\tilde{G}'$.

*Proof.* Let us assume first that n is even and let us denote by V and V' the representation spaces of ρ and ρ' respectively. According to 2.6, it is enough to prove that the categories $\mathcal{E}_G^V$ and $\mathcal{E}_G^{V'}$ are equivalent. To that effect, we apply the method in 1.10. If T is a space of spinorial representations of even dimension, the categories $\mathcal{E}_G^W$ and $\mathcal{E}_G^{T\oplus W}$ are equivalent. Moreover, we have V ⊕ T ≅ V' ⊕ T', with T = V ⊕ V ⊕ V' ⊕ V' and T' = V ⊕ V ⊕ V ⊕ V' of a Spin type. The case where n is odd is deduced from the even case if we replace V and V' by V ⊕ 1 and V' ⊕ 1 respectively.

**2.10.** The number of conjugacy classes of G or **Z**/2 ⋉ G which are decomposed in $\tilde{G}$ or **Z**/2 ⋉ $\tilde{G}$ is determined explicitly by theorem 2.4. If G is the symmetric group on n letters, these numbers may be computed in terms of suitable partitions of n. More precisely, in addition to the numbers $p_n$ and $i_n$ of the introduction, we define $j_n$ as the number of partitions of n of the type n = $\lambda_1$ +... + $\lambda_s$, 1 ≤ $\lambda_1$ ≤... ≤ $\lambda_s$, with odd $\lambda_i$'s. It is a classical fact that $j_n = p_n + i_n$. The reason is the following : the generating series of $p_n + i_n$ is the product $\prod_{m=1}^{\infty}(1 + x^m)$, while the generating series of $j_n$ is the product $\prod_{r=1}^{\infty}(1 - x^{2r-1})^{-1}$. Therefore, we only have to check the identity

$$\prod_{m=1}^{\infty}(1 + x^m)(1 - x^{2m-1}) = 1$$

which is quite obvious by a straightforward computation. From 2.5 and this discussion, we deduce the following statement which is an extension of the theorem 1.9 (see also theorem 3.10) :

**2.11. THEOREM.** *Let G be the symmetric group* $\mathfrak{S}_n$ *and* $V = \mathbf{R}^n$. *Then*

- $K(\mathcal{E}_G^V)$ *is a free group of rank* $p_n + 2i_n$ *and*
- $K(\mathcal{E}_G^{V \oplus 1})$ *is a free group of rank* $2p_n + i_n$

*Proof.* With the notations of the table at the end of this paper, where we look at all types of partitions of n, we have $p_n$ = B0 + B2 + D4 ; $i_n$ = D0 + D2 + B4 ; $j_n$ = A0 + B0 + C3 + D0 ( = $p_n + i_n$ according to 2.10). Note that the rows D (resp. B) are zero if n is even (resp. odd).

**2.12. Generalization.** Let A be a central simple algebra over a field k and let G be a finite group acting on A with an order invertible in A (which implies that the crossed product algebra $G \ltimes A$ is semi-simple). According to the Skolem-Noether theorem, the homomorphism $\Theta$ from A* to Aut(A) associating to $s \in A^*$ the inner automorphism $x \mapsto s.x.s^{-1}$ is surjective with kernel k*. With notations slightly different from the previous ones, we define the "linear Schur group" $\overline{G}$ by the following pull-back diagram

$$\begin{array}{ccc} \overline{G} & \longrightarrow & A^* \\ \downarrow & & \downarrow \\ G & \longrightarrow & \mathrm{Aut}(A) \end{array}$$

In this situation, a representation of $\overline{G}$ of linear type is defined as an ordinary representation of $\overline{G}$ such that $\rho(\lambda)$ is the multiplication by the scalar $\lambda$, for $\lambda \in k^* \subset \overline{G}$. With the same argument as in 2.6, we can show that the number of simple factors of the crossed product algebra $G \ltimes A$ is the number of irreducible representations of $\overline{G}$ of linear type. In our context, A is the Clifford algebra C(V), where V is of even dimension (over the complex field $\mathbf{C}$). In this situation the irreducible representations of $\overline{G}$ of linear type are in bijective correspondence with those of $\tilde{G} \subset \overline{G}$ (as in 2.5).

## 3. Equivariant K-theory of real projective spaces.

**3.1.** We may go further in the relation between Algebra and Topology by considering real projective bundles instead of Thom spaces. In so doing, we shall prove the result in 2.11 in a topological way. The following general theorem is valid in complex K-theory AND real K-theory : it is a natural extension of the results in [6] and [7].

**3.2. THEOREM.** *Let* V *be an equivariant real G-vector bundle with base* X *and* P(V) *its associated real projective bundle. The equivariant K-theory* $K_G^*(P(V))$ *is then naturally isomorphic to the Grothendieck group* $K^{*+1}(\phi)$, *where* $\phi$ *is the restriction functor*

$$\phi : \mathcal{E}_G^{V \oplus 1}(X) \longrightarrow \mathcal{E}_G^1(X) \ ,$$

*with "1" denoting as always the trivial vector bundle of rank* 1. *More generally, if* W *is an equivariant subbundle of* V, *the relative equivariant K-theory* $K_G^*(P(V), P(W))$ *is* *isomorphic to the Grothendieck group* $K^{*+1}(\psi)$, *where $\psi$ is the functor*

$$\psi : \mathcal{E}_G^{V\oplus 1}(X) \longrightarrow \mathcal{E}_G^{W\oplus 1}(X)$$

*Proof.* This theorem is proved in [6], p. 242, when G is the trivial group. Unfortunately, this proof does not generalize to the equivariant situation. We shall give here a more conceptual argument based on the following remark : the group $K_G^*(P(V))$ is naturally isomorphic to $K_{G \times \mathbb{Z}/2}^*(S(V))$, $S(V)$ being the sphere bundle with the usual antipodal action of $\mathbb{Z}/2$.

In order to simplify the notations, let us put $H = G \times \mathbb{Z}/2$. We have the cohomology exact sequence

$$\longrightarrow K_H^i(B(V), S(V)) \longrightarrow K_H^i(B(V)) \longrightarrow K_H^i(S(V)) \longrightarrow K_H^{i+1}(B(V), S(V)) \longrightarrow$$

In homotopical terms, this means that the group $K_H^i(S(V))$ may be identified with the $K^{i+1}$-group of the following diagram of Banach categories (cf. [6] p.192) :

$$\begin{array}{ccc} \mathcal{E}_H^{V\oplus 1}(X) & \longrightarrow & \mathcal{E}_H(X) \\ \downarrow & & \downarrow \\ \mathcal{E}_H^V(X) & \longrightarrow & 0 \end{array}$$

which is also the $K^{i+1}$-group of the obvious functor $\mathcal{E}_H^{V\oplus 1}(X) \longrightarrow \mathcal{E}_H(X) \times \mathcal{E}_H^V(X)$ [this general fact is true for any group H acting linearly on V]. Let us examine now the specific case where $H = G \times \mathbb{Z}/2$. Since $\mathbb{Z}/2$ acts on V by the involution $\varepsilon : v \mapsto -v$, the category $\mathcal{E}_H^V(X)$ is identified with $\mathcal{E}_G^{V\oplus 1}(X)$. On the other hand, one may describe the category $\mathcal{E}_H^{V\oplus 1}(X)$ by noticing that if the quadruple $(g, v, \varepsilon, \eta)$ describes the action of G V, $\mathbb{Z}/2$ and 1 respectively, we may associate to it bijectively the quadruple $(g, v, \varepsilon, \varepsilon\eta)$, with an involution $\varepsilon\eta$ commuting with the three other actions : this shows that the category $\mathcal{E}_H^{V\oplus 1}(X)$ is equivalent to the product $\mathcal{E}_G^{V\oplus 1}(X) \times \mathcal{E}_G^{V\oplus 1}(X)$. Therefore, the group $K_H^i(S(V))$ is isomorphic to the $K^{i+1}$-group of the functor

$$\mathcal{E}_G^{V\oplus 1}(X) \times \mathcal{E}_G^{V\oplus 1}(X) \longrightarrow \mathcal{E}_G^{V\oplus 1}(X) \times \mathcal{E}_G^1(X)$$

If we "eliminate" $\mathcal{E}_G^{V\oplus 1}(X)$, we recover the functor $\phi$ of the theorem.

The group $K_G^*(P(V), P(W))$ is determined in an analogous way : it is the $K^{*+1}$-group of the diagram

$$\begin{array}{ccc} \mathcal{E}_G^{V\oplus 1}(X) & \longrightarrow & \mathcal{E}_G^1(X) \\ \downarrow & & \downarrow \\ \mathcal{E}_G^{W\oplus 1}(X) & \longrightarrow & \mathcal{E}_G^1(X) \end{array}$$

which is also the $K^{*+1}$-group of the functor $\psi$ in the theorem.

**3.3.** Let us suppose now that X is reduced to a point[10]. We are going to determine the dimension of the vector space $K_G^*(P(V)) \otimes \mathbf{C}$, thanks to computations of the same nature as those of the first section. Making the terminology of 1.7 more precise, we say that a conjugacy class[11] $<g>$ is <u>strictly positive</u> (resp. <u>negative</u>) if $\text{Ker}(\rho(g) - 1) \neq 0$ (resp. $\text{Ker}(\rho(g) + 1) \neq 0$). Now let $\{g_i\}$ be the set of conjugacy classes of elements in G. If we put $X = P(V)$, the sub-vector space $X^{g_i}$ is the disjoint union of $P(V^{+g_i})$ and of $P(V^{-g_i})$, where

$$V^{+g_i} = \text{Ker}(\rho(g_i) - 1) \quad \text{and} \quad V^{-g_i} = \text{Ker}(\rho(g_i) + 1)$$

Therefore, we see that a conjugacy class which is counted strictly positive or negative[12] contributes one extra dimension in the group $K^0$.

In the same way, a conjugacy class which is even oriented, counted positively or negatively contributes for one dimension in the $K^1$-group (because the corresponding projective space is of odd dimension). From this discussion, we deduce the following theorem :

**3.4. THEOREM.** *Let $A_G$ (resp. $O_G$) be the number of conjugacy classes which are stricly positive or negative (resp. strictly positive or negative and even oriented). Then, the rank of the group $K_G^0(P(V))$ (resp. $K_G^1(P(V))$) is equal to $A_G$ (resp. $O_G$).*

**3.5. COROLLARY.** *Let $C_G$ be the number of conjugacy classes of G. Then the rank $R_G^{V\oplus 1}$ of the free group $K(\mathcal{E}_G^{V\oplus 1})$ satisfies the equation*

$$R_G^{V\oplus 1} - 2C_G = O_G - A_G$$

*Therefore, if V is odd dimensional (resp. even dimensional), $2C_G + O_G - A_G$ is the number of conjugacy classes of G (resp. $\mathbf{Z}/2 \times G$) which are decomposed in $\tilde{G}$ (resp. $\mathbf{Z}/2 \ltimes \tilde{G}$).*

*Proof.* It is a direct consequence of 2.8 et 3.3, since the category $\mathcal{E}_G^1$ is the product $\mathcal{E}_G \times \mathcal{E}_G$.

---

[10] We should notice that $K_G^*(P(V))$ is <u>not</u> free in general, in contrast with the group $K_G^*(V)$.

[11] A conjugacy class may be positive AND negative.

[12] A conjugacy class which is counted at the same time positive and negative contributes 2 dimensions in the associated K-group.

In the same way, we prove the following general statement :

**3.6. COROLLARY.** *Let $P_G$ (resp. $N_G$) be the number of strictly positive oriented conjugacy classes which are even (resp. odd). Then, according to 1.8, the rank $R_G^V$ of the free group $K(\mathcal{E}_G^V)$ satisfies the equation $R_G^{V\oplus 1} - R_G^V = P_G - N_G$, which may be also written*

$$R_G^V = 2C_G + O_G - A_G - P_G + N_G$$

*Therefore, if the dimension of V is even (resp. odd), $2C_G + O_G - A_G - P_G + N_G$ is the number of conjugacy classes of G (resp. $\mathbf{Z}/2 \times G$) which are decomposed in $\tilde{G}$ (resp. $\mathbf{Z}/2 \ltimes \tilde{G}$).*

**3.7.** Let us apply the previous considerations to the symmetric group $G = \mathfrak{S}_n$ acting naturally on $\mathbf{R}^n$. As we said before, the conjugacy class of g is associated to a partition $1 \leq \lambda_1 \leq ... \leq \lambda_s$ with $n = \lambda_1 + ... + \lambda_s$. In the table at the end of the paper, we describe all the types of partitions which are relevant and deduce from it and theorem 3.4 the following result[13]

**3.8. THEOREM.** *Let $G = \mathfrak{S}_n$ acting naturally on $V = \mathbf{R}^n$. Then the rank of $K_G^0(P(V))$ is $2P(n) - j_n$ and the rank of $K_G^1(P(V))$ is $p_n$, where P(n) is the total number of partitions of n, $p_n$ is defined in the introduction, and $j_n$ is the number of partitions of $n = \lambda_1 + ... + \lambda_s$ of type $1 \leq \lambda_1 \leq ... \leq \lambda_s$ with odd $\lambda_i$'s (note that $j_n = p_n + i_n$ according to 2.10).*

**3.9.** On the other hand, we have $P(n) = R_G^0$, the rank of $K(\mathcal{E}_G^0)$ which is the classical representation ring of G. In the same way, the rank $R_G^1$ of $K(\mathcal{E}_G^1) = K(\mathcal{E}_G^0 \times \mathcal{E}_G^0)$ is $2P(n)$. As a consequence of 3.2, we therefore get the following result which gives a more direct[14] topological proof of theorem 2.11 :

**3.10. THEOREM.** *The rank $R_G^{V\oplus 1}$ of the free group $K(\mathcal{E}_G^{V\oplus 1})$ is the solution of the equation*

$$R_G^{V\oplus 1} - 2P(n) = \mathrm{rank}(K_G^1(P(V))) - \mathrm{rank}(K_G^0(P(V))) = p_n - (2P(n) - j_n)$$

*and therefore*

$$R_G^{V\oplus 1} = p_n + j_n = 2p_n + i_n$$

*On the other hand, according to 1.1 and 1.3, the rank $R_G^V$ of the free group $K(\mathcal{E}_G^V)$ is the*

---

[13] We recall the notations of this table : $p_n$ = B0 + B2 + D4, $i_n$ = D0 + D2 + B4, $j_n$ = A0 + B0 + C3 + D0 (= $p_n + i_n$ according to 2.10).

[14] It is a more "direct" proof, since we do not distinguish between the even case from the odd.

*solution of the equation* $R_G^{V\oplus 1} - R_G^V = p_n - i_n$, *and therefore*

$$R_G^V = p_n + 2i_n$$

**3.11. Remark.** Theorem 3.10 has an equivalent formulation in terms of the K-theory of crossed product algebras. Let us denote by $C(\mathbf{R}^n)$ the Clifford algebra of $\mathbf{R}^n$ provided with a positive definite quadratic form. If $A_n$ is the crossed product algebra $\mathfrak{S}_n \ltimes C(\mathbf{R}^n)$ and $B_n$ the crossed product algebra $\mathfrak{S}_n \ltimes C(\mathbf{R}^{n+1})$, we have the following formulas

$$K(A_n) \cong \mathbf{Z}^{2p_n + i_n} \text{ et } K(B_n) \cong \mathbf{Z}^{p_n + 2i_n}$$

## 4. Relation with operations in complex K-theory. Multiplicative structures.

**4.1.** Let $R(\mathfrak{S}_n)$ be the representation ring of the symmetric group $G = \mathfrak{S}_n$. As Atiyah has shown in [2], every group homomorphism $R(G) \longrightarrow \mathbf{Z}$ gives rise to an operation in (complex) K-theory. It is defined by the composition $K(X) \to K_G(X) \cong K(X) \otimes R(G) \longrightarrow K(X)$, where the first map is $E \mapsto E^{\otimes n}$, the group G acting trivially on X and by permutation of the factors in $E^{\otimes n}$. Let us put $R(G)^* = \mathrm{Hom}(R(G), \mathbf{Z})$; the idea (due again to Atiyah) is to provide the direct sum $\oplus_n R(\mathfrak{S}_n)^*$ with an algebra structure (denoted by $\odot$) : it is induced by the obvious composition $R(\mathfrak{S}_m)^* \otimes R(\mathfrak{S}_n)^* \cong R(\mathfrak{S}_m \times \mathfrak{S}_n)^* \longrightarrow R(\mathfrak{S}_{m+n})^*$. As a key example, one might consider the element $\lambda^n$ of $R(\mathfrak{S}_n)^* = \mathrm{Hom}(R(\mathfrak{S}_n) \longrightarrow \mathbf{Z})$ defined by $\lambda^n(\pi) = 1$ if $\pi$ is the dimension 1 representation of $\mathfrak{S}_n$ associated to the signature and by $\lambda^n(\pi) = 0$ for the other irreducible representations.

**4.2. THEOREM** (Atiyah [2]). *The coalgebra $R(\mathfrak{S}_k)^*$ is a free $\mathbf{Z}$-module with basis the $\odot$-products $\lambda^{i_1} \odot \ldots \odot \lambda^{i_r}$ with $i_1 \leq i_2 \leq \ldots \leq i_r$ and $i_1 + \ldots + i_r = k$.*

**4.3**. We wish to make the coproduct $m : R(\mathfrak{S}_k)^* \longrightarrow R(\mathfrak{S}_k)^* \otimes R(\mathfrak{S}_k)^*$ more explicit. With obvious definitions, it is easy to see that $m(u \odot v) = m(u) \odot m(v)$. Theoreofore, it is enough to compute $m(\lambda^n)$, which amounts to computing $\lambda^n(xy)$ in terms of the $\lambda^i(x)$ and $\lambda^j(y)$ in the $\lambda$-ring $K(X)$. For instance, we have

$$\lambda^2(xy) = \lambda^2(x)(\lambda^1(y))^2 + (\lambda^1(x))^2 \lambda^2(x) - 2\lambda^2(x)\lambda^2(y)$$

and therefore

$$m(\lambda^2) = \lambda^2 \otimes (\lambda^1 \odot \lambda^1) + (\lambda^1 \odot \lambda^1) \otimes \lambda^2 - 2\lambda^2 \otimes \lambda^2$$

a formula which we may write in the condensed form

$$m(\lambda^2) = \lambda^2 \otimes (\lambda^1)^2 + (\lambda^1)^2 \otimes \lambda^2 - 2\lambda^2 \otimes \lambda^2$$

**4.4.** By tensoring all the previous modules by **Q**, we get much simpler formulas thanks to the Adams operations (cf. 4.6). More precisely, if $\lambda^1, \lambda^2, ..., \lambda^k$ are indeterminates, we denote by $\psi^i$ the Newton polynomial in the $\lambda^1, \lambda^2, ..., \lambda^k$ ($k \geq i$) which gives the fundamental symmetric function $S_i = \sum (x_r)^i$ in terms of the elementary symmetric functions $\sigma_1, ..., \sigma_i$. If $I = (i_1, ..., i_k)$ is such that $i_1 \leq i_2 \leq ... \leq i_k$, we denote by $\lambda^I$ (resp. $\psi^I$) the product $\lambda^{i_1} ... \lambda^{i_k}$ (resp. $\psi^{i_1} ... \psi^{i_k}$). We have formulas giving the $\psi^I$ in terms of the $\lambda^I$ and vice-versa:

$$\psi^J = \sum_I \theta^J_I \lambda^I$$

$$\lambda^I = \sum_J \omega^I_J \psi^J$$

Here $\theta^J_I$ (resp. $\omega^I_J$) is a "triangular" matrix with integral coefficients (resp. rational coefficients). If we put $|I| = i_1 + i_2 + ... + i_k$, the previous theorem has a "dual" formulation:

**4.5. THEOREM.** *The representation ring $R(\mathfrak{S}_k)$ is a free $\mathbf{Z}$-module with basis the "dual" symbols $\lambda_I$ with $|I| = k$. The product $\lambda_I . \lambda_J$ is given by the following formula*

$$\sum_{K, L} \theta^K_I \theta^K_J \omega^L_K \lambda_L$$

*In particular, the coefficient* $c^L_{IJ} = \sum_K \theta^K_I \theta^K_J \omega^L_K$ *is an integer. This is the coefficient $c^L_{IJ}$ in the universal formula for $\lambda$-rings when we write $\lambda^L(xy)$ in terms of products of exterior powers $\lambda^I(x)$ and $\lambda^J(y)$:*

$$\lambda^L(xy) = \sum_{I, J} c^L_{IJ} \lambda^I(x) \lambda^J(y)$$

**4.6. Remark.** It will be more convenient later on to consider the operations $\gamma_i$ defined by Grothendieck ([8] p. 253) instead of the $\lambda_i$'s. The structural equations $\gamma^L(xy) = \sum_{I, J} u^L_{IJ} \gamma^I(x) \gamma^J(y)$ in general $\lambda$-rings enable us to compute in the same way the multiplicative structure of $R(\mathfrak{S}_n)$ in terms of the basis given by the $\gamma^J$.

**4.7.** The same considerations apply to the $R(G)$-module $K^*_G(\mathbf{R}^n)$. More precisely, to each homomorphism $K^*_G(\mathbf{R}^n) \longrightarrow \mathbf{Z}$, where $G = \mathfrak{S}_n$, we may associate an operation in K-theory by the composition of the following maps

$$K^1(X) = K(X \times \mathbf{R}) \longrightarrow K_G(X \times \mathbf{R}^n) \cong \bigoplus_{r+s=0} K^r(X) \otimes K^s_G(\mathbf{R}^n) \longrightarrow K^*(X)$$

More precisely, an element of $\mathrm{Hom}(K_G^0(\mathbf{R}^n), \mathbf{Z})$ (resp. $\mathrm{Hom}(K_G^1(\mathbf{R}^n), \mathbf{Z})$) gives rise to an operation $K^1(X) \longrightarrow K^0(X)$ (resp. $K^1(X) \longrightarrow K^1(X)$). These operations are classical and determined by the K-theory of the infinite unitary group : we obtain a free $\mathbf{Z}$-module with basis the symbols $\gamma^S = \gamma^{s_1}\gamma^{s_2}...\gamma^{s_r}$ with $s_1 < s_2 < ... < s_r$ and $s_1 + ... + s_r = n$, the parity of r being equal to the degree of the operation. On the other hand, if we consider an element xy, where $x \in \mathrm{Hom}(R(G), \mathbf{Z})$ and $y \in \mathrm{Hom}(K_G^*(\mathbf{R}^n), \mathbf{Z})$, with the identification $\mathrm{Hom}(R(G) \otimes K_G^*(\mathbf{R}^n), \mathbf{Z}) \cong \mathrm{Hom}(R(G), \mathbf{Z}) \otimes \mathrm{Hom}(K_G^*(\mathbf{R}^n), \mathbf{Z})$, we have the same formula as above

$$\gamma^S(xy) = \sum_{I,J} u_{IJ}^S \gamma^I(x) \gamma^J(y)$$

However, we have to replace $u_{IJ}^S$ by 0 if J is a partition with non distincts integers (because the squares of elements in $K^1$ are equal to 0) and we therefore obtain the following theorem :

**4.8. THEOREM.** *The group $K_G^*(\mathbf{R}^n)$ is a free $\mathbf{Z}$-module with basis the "dual" symbols $\gamma_S$, where S is a partition of n into <u>distinct</u> integers and where the degree is determined by the parity of the number of terms in the partition. The R(G)-module structure of $K_G^*(\mathbf{R}^n)$ is given by the following relation*:

$$\gamma_I \cdot \gamma_S = \sum_{I,J} u_{IJ}^S \gamma_J$$

*where the integers $u_{IJ}^S$ are defined above and where S and J are partitions of n into <u>distinct</u> integers .*

**4.9. Remarks.** As in 4.4, we may take the tensor products of these R(G)-modules by $\mathbf{Q}$. That way, $K_G^*(\mathbf{R}^n) \otimes \mathbf{Q}$ is identified with the sum of the ideals in $R(G) \otimes \mathbf{Q}$ given by the idempotents $\psi_S$, S running through the set of partitions of n into <u>distinct</u> integers. If G is an arbitrary finite group, the same argument may be applied to the R(G)-module $K_G^*(V) \otimes \mathbf{C}$, thanks to the equivariant Chern character [4] [10]. It gives an isomorphism between the following R(G)-modules

$$K_G^*(V) \otimes \mathbf{C} \cong \bigoplus_{g_i \in <G>} H_c^*(V^{g_i}/C_{g_i} ; \mathbf{C})$$

Note that the R(G)-module structure of the second factor comes from the fact that $R(G) \otimes \mathbf{C}$ is isomorphic to the ring of complex valued functions on the set of conjugacy classes of G.

Equipe Topologie et Géométrie Algébriques
Institut de Mathématiques de Jussieu
Université Paris 7-Denis Diderot - Case 7012
2, Place Jussieu, 75251 PARIS cedex 05 - FRANCE
e.mail and Web : karoubi@math.jussieu.fr ; http://www.math.jussieu.fr/~karoubi/